\documentclass{amsart}
\usepackage{amsfonts,amssymb,enumerate}

\newcommand{\N}{\ensuremath{\mathbb{N}}}
\newcommand{\e}{\varepsilon}
\newcommand{\clo}{\overline}
\newcommand{\norm}[1]{{\Vert #1 \Vert}}
\newcommand{\dens}{{\rm dens}}

\newtheorem{theorem}{Theorem}
\newtheorem{lemma}[theorem]{Lemma}
\newtheorem{corollary}[theorem]{Corollary}
\newtheorem{proposition}[theorem]{Proposition}

\begin{document}

\title{Connected, not separably connected complete metric spaces}
\author{T.~Banakh, M.~Vovk, M.~R.~W\'ojcik}

\address{Taras Banakh: Ivan Franko Lviv National University, Lviv, Ukraine, and\newline
Unwersytet Humanistyczno-Przyrodniczy im. Jana Kochanowskiego, Kielce, Poland}
\email{tbanakh@yahoo.com}
\address{M.~Vovk: National University ``Lvivska Politechnika", Lviv, Ukraine}
\address{Micha{\l} Ryszard W\'ojcik}
\email{michal.ryszard.wojcik@gmail.com}

\begin{abstract}
In a separably connected space any two points
are contained in a separable connected subset.
We show a mechanism that takes a connected bounded metric space
and produces a complete connected metric space whose
separablewise components form a quotient space
isometric to the original space.
We repeatedly apply this mechanism to construct,
as an inverse limit, a complete connected
metric space whose each separable subset
is zero-dimensional.
\end{abstract}

\maketitle

A topological space is {\bf separably connected} iff
any two of its points are contained in a separable connected subset.
A {\bf separablewise component} is the union of all separable connected
subsets containing a given point.
Since the closure of a separable connected set is again separable and connected,
it follows that separablewise components are closed.
A connected space whose all connected separable subsets are singletons
will be called {\bf nonseparably connected}.

The first example of a nonseparably connected metric space was given by Pol
in 1975, \cite{POL}. Another example was given by Simon in 2001, \cite{SIMON}.
In 2003, Aron and Maestre constructed  a connected,
not separably connected metric space that contains many arcs, \cite{ARON}.
In 2008, Morayne and W\'ojcik obtained a nonseparably
connected metric space as a graph of a function from the real line
satisfying Cauchy's equation and thus forming a topological group,
\cite{PhD} or \cite{MW}.
None of these spaces are completely metrizable.

Our first result is that
for any connected metric space $X$ we produce a complete connected metric space
whose separablewise components form a quotient space homeomorphic to $X$,
Corollary \ref{abundance}.
(This yields a lot of topologically different
completely metrizable connected, not separably connected spaces.)

Our second result is the first known example
of a nonseparably connected metric space that is complete,
Corollary \ref{firstever}.
In fact, our space satisfies a stronger property.
Namely, the metric is economical in the following sense.

\subsection*{Economically metrizable spaces}
Given a metric space $(X,d)$, we say that the metric $d$ is {\em economical} iff
$$|d(A\times A)|\le \dens(A)$$
for any infinite subset $A\subset X$,
where $$d(A\times A)=\{d(a,b)\colon a,b\in A\},$$
$$\dens(A)=\min\{|D|\colon D\subset A\subset\clo D\}.$$
\begin{proposition}\label{prop1}
Every separable subset of an economically metrizable space
is zero-dimensional.
\end{proposition}
\begin{proposition}\label{prop2}
Every connected economically metrizable space is nonseparably connected.
\end{proposition}

\begin{proposition}
The Cantor set is economically metrizable.
\end{proposition}
\begin{proof}
Let $\N=\{1,2,3,\ldots\}$.
Consider the following metric $d$ on the space $P(\N)$ of all subsets of $\N$,
$$d(A,B)=
\begin{cases}0&\text{ if }A=B,\\
\cfrac{1}{\min((A\setminus B)\cup(B\setminus A))}&\text{ if }A\not=B.
\end{cases}$$
Now, the economically metrizable space $(P(\N),d)$ is homeomorphic
to the ternary Cantor set with the euclidean metric,
which is not economical.
\end{proof}

\begin{proposition}
Every metrizable space containing a copy of the Cantor set
admits a metric that is not economical.
\end{proposition}
\begin{proof}
Recall that given a metrizable space $X$ and a closed subset $M\subset X$,
every admissable metric on $M$ can be extended to an admissable metric on $X$,
Engelking 4.5.21(c).
So if $X$ contains a copy of the Cantor set, $M\subset X$,
we may choose a metric for $M$ that is not economical
and extend it to the whole $X$.
\end{proof}

\section{Preliminaries: connected inverse limits}

In this section we prepare the tools that we use
to show that the spaces we construct are connected.
Lemma \ref{openineachfiber} and Lemma \ref{ZconnectedifXconnected}
are used in Theorem \ref{supermachina}
to obtain connected, not separably connected spaces.
The remaining lemmas and Theorem \ref{Puzio}
are used in Theorem \ref{economical} 
to show that our inverse limit is connected.

We need to make a number of simple observations ---
Lemmas \ref{supermonotone}, \ref{stillmonotone}, \ref{herquocompo}
--- before we are ready to apply E. Puzio's Theorem 11
from \cite{Puzio}.
\\\\
Let $X$ and $Y$ be topological spaces.
Then a function $f\colon X\to Y$ is
\begin{enumerate}[$\bullet$]
\item {\em monotone} iff
$f^{-1}(y)$ is connected for each $y\in Y$,
\item {\em open at} $x\in X$
iff $f(x)\in Int(f(U))$ for every neighborhood $U$ of $x$,
\item {\em hereditarily quotient}
iff $$f^{-1}(y)\subset U\implies y\in Int(f(U))$$
for every open set $U\subset X$ and every point $y\in Y$.
\end{enumerate}

\begin{lemma}\label{openineachfiber}
If every fiber of $f\colon X\to Y$ contains
a point at which $f$ is open,
then $f$ is hereditarily quotient.
\end{lemma}
\begin{proof}
This is immediate from the definitions.
\end{proof}

\begin{lemma}[Engelking 6.1.I]
\label{ZconnectedifXconnected}
Suppose that $f\colon X\to Y$ is monotone and
hereditarily quotient.
Then X is connected if f(X) is connected.
\end{lemma}
\begin{proof}
Let $U$ be a clopen subset of $X$.
To show that $f(U)$ is open take any $y\in f(U)$.
Since the fiber $f^{-1}(y)$ is connected
and intersects the clopen set $U$,
it follows that $f^{-1}(y)\subset U$.
Since $f$ is hereditarily quotient, $y\in Int(f(U))$,
showing that $f(U)$ is open.
Since $X\setminus U$ is also clopen,
it follows by analogy that $f(X\setminus U)$ is open, too.
Now, the sets $f(U)$, $f(X\setminus U)$ are disjoint,
because the clopen sets $U$, $X\setminus U$ contain whole fibers.
The connected space $f(X)=f(U)\cup f(X\setminus U)$
is a union of two disjoint open subsets.
Thus $U=X$ or $U=\emptyset$, showing that $X$ is connected.
\end{proof}

\begin{lemma}\label{supermonotone}
If $f\colon X\to Y$ is monotone and hereditarily quotient
then $f^{-1}(E)$ is connected whenever $E\subset Y$ is connected.
\end{lemma}
\begin{proof}
A straightforward subspace topology argument shows that the restriction
$$f|_{f^{-1}(E)}\colon f^{-1}(E)\to E$$ is hereditarily quotient.
It is also evidently monotone.
Therefore, by Lemma~\ref{ZconnectedifXconnected},
$f^{-1}(E)$ is connected if $E$ is connected.
\end{proof}

\begin{lemma}\label{stillmonotone}
If $f\colon X\to Y$ and $g\colon Y\to Z$ are monotone
and $f$ is hereditarily quotient,
then $g\circ f$ is monotone.
\end{lemma}
\begin{proof}
Since $g^{-1}(z)$ is connected
and $f$ is monotone and hereditarily quotient,
by Lemma~\ref{supermonotone},
$(g\circ f)^{-1}(z)=f^{-1}(g^{-1}(z))$ is connected.
\end{proof}

\begin{lemma}\label{herquocompo}
If $f\colon X\to Y$ and $g\colon Y\to Z$ are hereditarily quotient
then their composition $g\circ f\colon X\to Z$ is hereditarily quotient.
\end{lemma}
\begin{proof}
Take any open set $U\subset X$ and any point $z\in Z$
such that $(g\circ f)^{-1}(z)\subset U$. Then
$$(g\circ f)^{-1}(z)=f^{-1}(g^{-1}(z)))=\bigcup_{y\in g^{-1}(z)}f^{-1}(y)\ \subset\ U.$$
Since $f$ is hereditarily quotient, we have $$y\in Int(f(U))$$
for each $y\in g^{-1}(z)$. Thus $g^{-1}(z)\subset Int(f(U))$.
Now, since $g$ is hereditarily quotient,
$$z\in Int(g(Int(f(U))))\subset Int((g\circ f)(U)).$$
\end{proof}

\begin{theorem}[E. Puzio, 1972]\label{Puzio}
Let $X_1,X_2,\ldots$ be a sequence of connected spaces.
Suppose that each function $f_n\colon X_{n+1}\to X_n$
is a continuous monotone hereditarily quotient surjection.
Then the inverse limit of this system
$$X=\big\{(x_n)_{n=1}^\infty\in\prod_{n=1}^\infty X_n\colon
(\forall n\in\mathbb N)\big(x_n=f(x_{n+1})\big)\big\}$$
is connected.
\end{theorem}
\begin{proof}
Theorem 11 in \cite{Puzio}
requires the additional assumption that the functions
$$f_n\circ f_{n+1}\circ\ldots\circ f_m\colon X_{m+1}\to X_n$$
are monotone hereditarily quotient surjections
for all $n<m\in\mathbb N$.
But this follows from our assumptions
thanks to Lemma~\ref{stillmonotone} and Lemma~\ref{herquocompo}.
So $X$ is connected.
\end{proof}

The little example $[0,1)\cup\{2\}\rightarrow[0,1]$
shows that the assumption that $f$ is hereditarily quotient
is needed in all these lemmas.

\section{Preliminaries: first countable spaces}

Our mechanism (Theorem \ref{supermachina})
that returns an appropriate complete metric space
for a given metric space actually works for any first countable space,
so we decided to write it more generally
at the cost of introducing some technical details for
dealing with first countable spaces.

\begin{proposition}\label{firstcountabled}
A topological space $X$ is first countable if and only if
to each pair of points $x,u\in X$ we may assign a number
$d_x(u)\in[0,1]$ so that
$$\inf\{d_x(u)\colon u\in E\}=0\iff x\in\overline E$$
for all $x\in X$ and $E\subset X$.
\end{proposition}
\begin{proof}
If $X$ is first countable, there are open sets
$\{U_n(x)\colon x\in X,n\in\mathbb N\}$ such that
$x\in U_{n+1}\subset U_n(x)$
and the sequence $U_n(x)$ is a local basis at $x$.
Then $$d_x(u)=\inf\{1/n\colon u\in U_n(x)\}$$
is easily seen to be as desired.

On the other hand, if $X$ admits such a function $d_x(u)$,
let $$B(x,r)=\{u\in X\colon d_x(u)<r\}$$
for all $x\in X,r>0$. Then
$$\inf\big\{d_x(u)\colon u\in X\setminus B(x,r)\big\}\ge r>0$$
and consequently, $x\not\in\overline{X\setminus B(x,r)}$.
In other words, $x\in Int(B(x,r))$.

If $G$ is an open set containing $x$,
then $$\inf\{d_x(u)\colon u\in X\setminus G\}\ge r>0,$$
for some $r>0$ and consequently $B(x,r)\subset G$.
So, $U_n(x)=Int(B(x,1/n))$ are the open sets which show
$X$ to be first countable.
\end{proof}

\section{Complete connected, not separably connected spaces}

\subsection*{Locally constant functions}
Let $X$ be a topological space. We say that a function
$f\colon X\to Y$ is {\em locally constant} iff
each point $x\in X$ has a neighborhood $U_x$ such that
$f|U_x$ is constant.
Naturally, the cardinality of the image $f(X)$
cannot exceed the density of $X$,
$|f(X)|\le \dens(X)=\min\{|D|\colon\clo D=X\}$.

\subsection*{Metrically discrete subsets}
Let $(X,d)$ be a metric space. We say that $A\subset X$ is
a {\em metrically discrete} subset of $X$
iff there is an $\e>0$ such that
$d(a,b)\ge\e$ for any distinct $a,b\in A$.
Naturally, the cardinality of a metrically discrete subset
cannot exceed the density of the metric space $X$,
$|A|\le \dens(X)$.

\begin{lemma}\label{spreadlemma}
If $X$ is a metric space and $f\colon X\to Y$ is locally constant
except on a metrically discrete subset,
then $|f(X)|\le \dens(X)$.
\end{lemma}
\begin{proof}
Let $K,D$ be disjoint subsets of $X$ such that
$f|K$ is locally constant and $D$ is metrically discrete.
Then $$|f(X)|\le|f(K)|+|f(D)|\le \dens(X)+|D|\le \dens(X)+\dens(X).$$
\end{proof}

\subsection*{Functionally Hausdorff spaces}
We say that a topological space $X$ is {\em functionally Hausdorff}
iff for any two distinct points $a,b\in X$ there is a continuous
function $f\colon X\to\mathbb R$ with $f(a)\not=f(b)$.
Each connected subset $E$ of such a space is either a singleton or contains
a set of cardinality $\mathfrak c$, as the following easy argument shows.
Let $a,b\in E$ with $f(a)<f(b)$.
Then the connected set $f(E)$ must contain the interval $[f(a),f(b)]$.

\begin{lemma}\label{NotSeparablyConnected}
Let $X$ be a metric space.
Let $Y$ be a functionally Hausdorff space.
Let $f\colon X\to Y$ be a Darboux function
that is locally constant except on a metrically discrete subset.
Then $f$ is constant on every connected separable subset of $X$.
Moreover, if $f$ is not constant, then $X$ is not separably connected.
\end{lemma}
\begin{proof}
Let $E\subset X$ be a connected separable subset.
Since $f$ is Darboux and $E$ is connected,
the set $f(E)$ is a connected subset of a functionally Hausdorff space.
Thus $f(E)$ is either a singleton or $|f(E)|\ge\mathfrak c$.
By Lemma~\ref{spreadlemma}, $|f(E)|\le \dens(E)=\aleph_0$,
so the set $f(E)$ is countable and thus a singleton.
\end{proof}

To argue that the spaces constructed in Theorem \ref{supermachina}
are completely metrizable, we obtain them as closed subsets of
a certain canonical space (which we call the cobweb)
that is easily seen to be complete.

\subsection*{The hedgehog}
Let $\kappa$ be a cardinal number.
A hedgehog with $\kappa$ spikes, each of length $\e>0$,
is the space $H=\{(0,0)\}\cup(\kappa\times(0,\e])$
equipped with the metric $\rho$ given by
\[\rho((x,t),(u,s))=\begin{cases}
|t-s|&\text{if}\ \ x=u,\\
\ t+s&\text{if}\ \ x\not=u.
\end{cases}\]
It is easy to see that it is a complete
one-dimensional arcwise connected metric space.

\subsection*{The cobweb}
Let $\e>0$ be a fixed number.
Let $V$ be a subset of a normed space with $\norm{u-v}=\e$
for all distinct $u,v\in V$.
Let $W=\bigcup\{[u,v]\colon u,v\in V\}$,
where $[u,v]$ denotes the line segment between vectors $u,v$.
Let $\rho(a,b)$ be the infimum over finite sums
$$\sum_{i=1}^n\norm{x_i-x_{i-1}}$$
where $a=x_0$ and $x_n=b$ and
$$(\forall i)(\exists u,v\in V)\big(\{x_i,x_{i-1}\}\subset[u,v]\big).$$
Then $\rho$ is clearly a metric on $W$,
the distance being measured along the {\em threads},
so that at each vortex the space $(W,\rho)$ is locally isometric
to a hedgehog with appropriately shortened spikes,
and at the remaining points it is
locally isometric to appropriately short euclidean intervals.
The metric space $(W,\rho)$ will be called
{\em the cobweb spun over $V$}. It is clearly zero-dimensional.
{\bf The cobweb is complete.}
Indeed, if $x_n$ is a Cauchy sequence,
we have a $k\in\mathbb N$ such that
$\rho(x_n,x_k)\le\e/4$ for all $n\ge k$.
It is evident that there is a vortex $u\in V$ with $\rho(x_k,u)\le\e/2$.
Then $\rho(x_n,u)\le\rho(x_n,x_k)+\rho(x_k,u)\le 3\e/4$ for all $n\ge k$.
So the Cauchy subsequence $x_k,x_{k+1},\ldots$ is contained
in the hedgehog with vortex $u$ and spikes of length $3\e/4$,
so it converges because the hedgehog is complete.
\\\\
The following theorem is written so as to reveal
all the interesting properties of the spaces constructed
and to allow first countable functionally Hausdorff spaces
in the place of metric spaces.
However, some of the conclusions are stronger for metric spaces
and these are summarized in Corollary \ref{abundance}.

\begin{theorem}\label{supermachina}
Let $X$ be a first countable functionally Hausdorff space.
Then there exists a complete one-dimensional metric space $Z$
and a continuous monotone hereditarily quotient surjection $f\colon Z\to X$ such that
\begin{enumerate}[(1)]
\item each fiber of $f$ is homeomorphic to a hedgehog
with $|X|$-many spikes,\\
thus $|Z|=\mathfrak c|X|$
\item each fiber $Z_a=f^{-1}(a)$ has exactly one point $a^*\in Z_a$
at which $f$ is open;\\these points form a metrically discrete subset
\item $Z$ is connected $\iff X$ is connected
\item $Z_a\setminus\{a^*\}$ is open in $Z$ for each $a\in X$,
thus $\dens(Z)\ge|X|$\\and $f$ is locally constant except at those discretely spaced points $a^*$
\item $|f(A)|\le \dens(A)$ for any subset $A\subset Z$
\item each connected separable subset of $Z$ lies in one of the fibers
\item $Z$ is not separably connected
\item its arcwise components coincide with the fibers
\item the space $Z|_X=\{Z_a\colon a\in X\}$ with $E\subset Z|_X$ declared open
iff $\bigcup E$ is open in $Z$ is homeomorphic to $X$ via
$Z|_X\ni Z_a\mapsto a\in X$.
\item $Z$ is locally connected at $z\in Z\iff z\notin\{a^*\colon a\in X\}$
\end{enumerate}
\end{theorem}
\begin{proof}
Since $X$ is first countable, let $d$ be as in Proposition~\ref{firstcountabled}.
Let $\overline{X}=X\times\{0\}$ and let $X^*=X\times\{1\}$.
Let $\Omega=\overline X\cup X^*$ be considered with the counting measure.
Let $Y$ be the Hilbert space $L^2(\Omega)$, that is
$$Y=\big\{h\colon\Omega\to\mathbb R\colon
\sum_{\omega\in\Omega}|h(\omega)|^2<\infty\big\}.$$
For each $x\in X$, let $\overline{x}\in Y$ be the characteristic
function of the singleton $\{(x,0)\}$, and let $x^*\in Y$ be the
characteristic function of the singleton $\{(x,1)\}$.
Let $V=\{\overline{x}\colon x\in X\}\cup\{x^*\colon x\in X\}$.
Then $\norm{u-v}=\sqrt 2$ for distinct $u,v\in V$.

Let $[u,v]$ denote the line segment joining two vectors $u,v\in Y$.
For any distinct $a,b\in X$, let $\overline{a}_b$ denote the point
lying in $[\overline{a},b^*]$ with $\norm{\overline{a}_b-b^*}=d_b(a)\le 1$.
Since $X$ is Hausdorff, $d_b(a)>0$, and so $\overline a_b\not=b^*$.
Let $$Z_a=[\overline{a},a^*]\cup
\bigcup\{[\overline{a},\overline{a}_b]\colon b\in X,b\not=a\}$$
for each $a\in X$. Let $Z=\bigcup\{Z_a\colon a\in X\}$.
Let $Z$ be equipped with the induced metric from the cobweb spun over $V$.

Notice that $Z$ is obtained by taking away selected open intervals
from some of the threads of the cobweb spun over $V$.
So $Z$ is a closed subset of this cobweb.
Thus $Z$ is a complete zero-dimensional metric space.

Let $f\colon Z\to X$ be given by $f(Z_a)=\{a\}$ for all $a\in X$.

Clearly, $f^{-1}(a)$ is connected for each $a\in X$, so $f$ is monotone.

We claim that $f$ is open at each $a^*\in Z$ because
$$f(a^*)=a\in Int(B(a,r))\subset f(B(a^*,r))$$ for all $a\in X,r>0$.
Indeed, for any $b\in B(a,r)$, we have
$\norm{\overline b_a-a^*}=d_a(b)<r$,
and so $b=f(\overline b_a)\in f(B(a^*,r))$.

By Lemma~\ref{openineachfiber}, $f$ is hereditarily quotient.

By Lemma~\ref{ZconnectedifXconnected}, $Z$ is connected if $X$ is connected.
Later we show that $f$ is continuous, so $X$ is connected if $Z$ is connected.

For each $a\in X$, the set $Z_a\setminus\{a^*\}$ is open in $Z$.
Thus $f$ is locally constant except on the metrically discrete set
$\{a^*\colon a\in X\}\subset V$. Therefore, for any set $A\subset Z$,
the restriction $f|A$ is locally constant except possibly
on a metrically discrete set, and by Lemma~\ref{spreadlemma}, $|f(A)|\le \dens(A)$.

To show that $f$ is continuous, it is sufficient to
analyze those points where $f$ is not locally constant.
Take any $a\in X$ and any $s\in(0,\sqrt 2)$.
Then, by Proposition~\ref{firstcountabled}, $f(a^*)\in Int(B(a,s))$,
and there is an $r\in(0,s)$ with $B(a,r)\subset Int(B(a,s))$.
We claim that $f$ is continuous at $a^*$ because
$$f(B(a^*,r))\subset B(a,r)\subset Int(B(f(a^*),s).$$
Indeed, take any $z\in B(a^*,r)$.
If $z\in Z_a$, then $f(z)=a\in B(a,r)$.
If $z\in Z_b$ with $b\not=a$, then $z\in[\overline b,\overline b_a]$
because $r$ is sufficiently small. Consequently,
$\norm{z-\overline b_a}+\norm{\overline b_a-a^*}=\norm{z-a^*}<r$.
Thus $d_a(b)=\norm{\overline b_a-a^*}<r$, and so $f(z)=b\in B(a,r)$.
Hence $f$ is continuous.

Now, $f$ is a Darboux function that is locally constant
except on a discrete subset.
Since $X$ is functionally Hausdorff,
by Lemma~\ref{NotSeparablyConnected},
we conclude that each connected separable subset of $Z$
lies in one of the fibers.
This means that the fibers coincide with
the separablewise components of $Z$,
each of which is homeomorphic to a hedgehog
and thus arcwise connected.
\end{proof}

\begin{corollary}\label{abundance}
Let $X$ be a metric space bounded by one.
Then there exists a complete metric space $Z$
(connected if and only if $X$ is connected)
and a Lipschitz monotone hereditarily quotient surjection
$f\colon Z\to X$ whose fibers coincide with the
separablewise components of $Z$ and form a quotient space
isometric to $X$.
\end{corollary}
\begin{proof}
Notice that in the course of the proof of Theorem~\ref{supermachina} we have
$$\min\{\rho(x,u)\colon x\in Z_a,u\in Z_b\}=\min\{d_a(b),d_b(a)\}.$$
Then the metric $z(Z_a,Z_b)=d(a,b)$ is an isometry because
$$f(B(a^*,r))=B(a,r)$$
for all $a\in X,r>0$.
We also have
$$d(f(x),f(u))\le\rho(x,u)$$
for all $x,u\in Z$.
\end{proof}

In the next section we obtain a complete nonseparably connected space
as an inverse limit of a sequence of spaces generated by
repeatedly applying Theorem \ref{supermachina}.
The following Corollary \ref{uniquesequence} lists
only those properties that are essential for that purpose.

\begin{corollary}\label{uniquesequence}
To every metric space $X$ there corresponds a uniquely determined
complete connected metric space $Z$ and a continuous monotone
hereditarily quotient surjection $f\colon Z\to X$ such that
$|f(A)|\le \dens(A)$ for any subset $A\subset Z$.
\end{corollary}

Recall that a connected, locally connected complete metric space
must be arcwise connected, Engelking 6.3.11.
The space $Z$ obtained in Theorem \ref{supermachina}
is not locally connected, although it is locally connected except on
a metrically discrete subset. This illustrates how important it is
to assume that the space is locally connected at each point
if we want to conclude that it is arcwise connected.

\section{A complete nonseparably connected space}

We make an appropriate choice of the product metric
for our inverse limit to ensure that it is 
economically and completely metrizable at the same time.

\begin{lemma}\label{dAxA}
Let $(X_n,d_n)$ be a sequence of uniformly bounded metric spaces.
Let $X=\prod_{n\in\N}X_n$ and let $\pi_n\colon X\to X_n$ be given by
$\pi_n(x_1,\ldots,x_n,\ldots)=x_n$.
Then $$d(x,u)=\max\Big\{\cfrac{d_n(\pi_n(x),\pi_n(u))}{2^{n}}\colon n\in\N\Big\}$$
for each $x,u\in X$ defines a product metric that is complete
if all factor metrics are complete.
Moreover, for any infinite subset $A\subset X$ we have
$$|d(A\times A)|\le\sup_{n\in\N}|\pi_n(A)|.$$
\end{lemma}
\begin{proof}
For any $a,b\in X$ we have
$$d(a,b)\in\bigcup\{2^{-n}d_n(\pi_n(a),\pi_n(b))\colon n\in\N\}.$$
Therefore, if $A\subset X$ is infinite, we have
$$|d(A\times A)|\le\sum_{n\in\N}|d_n(\pi_n(A)\times\pi_n(A))|\le
\sum_{n\in\N}|\pi_n(A)|^2\le\sup_{n\in\N}|\pi_n(A)|.$$
\end{proof}

\begin{theorem}\label{economical}
There exists a complete connected economical metric space.
\end{theorem}
\begin{proof}
Let $X_1$ be a connected complete metric space bounded by one.
By Corollary~\ref{uniquesequence}, we have a sequence
$(X_n,d_n)$ of complete connected metric spaces bounded by one
and continuous monotone hereditarily quotient surjections
$f_n\colon X_{n+1}\to X_n$ such that
$|f_n(A)|\le \dens(A)$ for any subset $A\subset X_{n+1}$.

Let $X$ be the inverse limit of this system.
By Theorem~\ref{Puzio}, $X$ is connected.

Equipping the product $\prod_{n=1}^\infty X_n$
with the product metric from Lemma~\ref{dAxA}
induces a complete metric $d$ on $X$ that satisfies
$|d(A\times A)|\le\sup_{n\in\N}|\pi_n(A)|$
for any infinite subset $A\subset X$.
It remains to show that $(X,d)$ is economical.
Let $A\subset X$ be any infinite subset.
We have $$|\pi_n(A)|=|f_n(\pi_{n+1}(A))|\le \dens(\pi_{n+1}(A))\le \dens(A)$$
for all $n\in\mathbb N$. Therefore,
$$|d(A\times A)|\le\sup_{n\in\N}|\pi_n(A)|\le \dens(A).$$
\end{proof}

\begin{corollary}\label{firstever}
There exists a complete connected metric space whose
each separable subset iz zero-dimensional.
In particular,
there exists a complete nonseparably connected metric space.
\end{corollary}
\begin{proof}
By Theorem \ref{economical} and Proposition \ref{prop1}.
\end{proof}

\subsection*{Connected punctiform spaces}
Recall that a topological space is {\em punctiform}
if all of its connected compact subsets are singletons.
For example, any Bernstein subset of the euclidean plane is
a connected punctiform metric space.
A separable complete connected punctiform space was constructed
by Kuratowski and Sierpi{\'n}ski in 1922, \cite{connexes_punctiformes}.
A variation of their idea was presented in \cite{PhD}.
Our complete nonseparably connected space
is a new example of a connected punctiform space.

\section{Appendix: a non-constant continuous locally extremal function}

Alessandro Fedeli and Attilio Le Donne constructed a connected metric space $X$
and a non-constant continuous function $f\colon X\to[0,1]$ that
has a local maximum or a local minimum at every point,
without claiming that the domain $X$ is complete, \cite{DF1} and \cite{DF2}.
By modifying the proof of our Theorem \ref{supermachina}, we make a similar construction
that is less technically burdensome, and the domain of the function
is easily seen to be complete.

\begin{theorem}\label{nonconstantcudo}
There is a complete connected metric space $(X,\rho)$
and a Lipschitz monotone hereditarily quotient surjection
$f\colon X\to(0,1)$ that has a local extremum at every point.
\end{theorem}
\begin{proof}
Let $\overline{a},a^\uparrow,a^\downarrow$ for $a\in(0,1)$
be identified as distinct points forming a discrete set $V$
in some nonseparable normed vector space, e.g.
$\norm{u-v}=1$ for $u\not=v$.
\\For $a\in(0,1)$ and $b\in(a,1)$,
let $b^\uparrow_a$ be the point lying in $[\overline a,b^\downarrow]$
with $\norm{b^\uparrow_a-b_b^\downarrow}=b-a$.
\\For $a\in(0,1)$ and $b\in(0,a)$,
let $b^\downarrow_a$ be the point lying in $[\overline a,b^\uparrow]$
with $\norm{b^\downarrow_a-b_b^\uparrow}=a-b$.
Let
$$H_a=[\overline{a},a^\uparrow]\cup[\overline{a},a^\downarrow]\cup
\big\{[\overline{a},b^\uparrow_a]\colon b\in(a,1)\big\}\cup
\big\{[\overline{a},b_a^\downarrow]\colon b\in(0,a)\big\}$$
for each $a\in(0,1)$. Let $X=\bigcup\{H_a\colon a\in(0,1)\}$.
Notice that $X$ is a closed subset of the cobweb spun over $V$.
Equipped with the metric $\rho$ induced from the cobweb,
$X$ is a complete metric space.

Let $f\colon X\to(0,1)$ be given by $f(H_a)=\{a\}$ for each $a\in(0,1)$.
Notice that $|f(x)-f(u)|\le\rho(x,u)$ for all $x,u\in X$,
so $f$ is continuous.
For each $a\in(0,1)$,
$f$ has a local minimum at $a^\uparrow$
and a local maximum at $a^\downarrow$.
It is locally constant at every other point.

Now, $f$ is hereditarily quotient because
$$f(B(a^\uparrow,r))=[a,a+r)$$
$$f(B(a^\downarrow,r))=(a-r,a]$$
for all $a\in(0,1)$ and all $r\in(0,a)\cap(0,1-a)$.
Since $f$ is a monotone hereditarily quotient surjection onto $(0,1)$,
by Lemma~\ref{ZconnectedifXconnected}, $X$ is connected.
\end{proof}

\subsection*{Further developments}
We are preparing a sequel paper in which we construct
a complete nonseparably connected space as a metric group, \cite{BW}.

Moreover, we plan to devote a separate paper to the construction
of the cobweb $\circledast(X)$ over a first countable space $X$,
having the properties listed in Theorem \ref{supermachina}.
We feel that this operation of obtaining the space $Z$
out of $X$ is interesting in itself
and needs to be investigated more closely, outside the context
of constructing connected, not separably connected spaces.
We have a number of different ways of describing this cobweb operation,
each having its advantages and disadvantages, and we feel that
this topic deserves separate treatment.

\subsection*{Acknowledgements}
We would like to thank Pawe{\l} Krupski for his topological seminar
at which we had a chance to present our work and improve it
greatly along the lines suggested by Krzysztof Omiljanowski.
\\\\
We rely on \cite{Puzio} in our Theorem \ref{economical}
to argue that the inverse limit is connected. Besides that,
the proofs of our main results are self-contained.

\end{document}